\documentclass[12pt]{article}
\usepackage{amsmath, amssymb, amsthm, enumerate}

\usepackage[utf8]{inputenc}
\usepackage{graphicx,color}

\theoremstyle{plain}
\newtheorem{theorem}{Theorem}[section]
\newtheorem{lemma}[theorem]{Lemma}
\newtheorem{proposition}[theorem]{Proposition}

\theoremstyle{definition}

\newtheorem{remark}[theorem]{Remark}

\def\rr{{\mathbb R}}

\def\sik{{\rr}^2}

\def\cc{{\mathbb C}}

\def\zz{{\mathbb Z}}

\def\gg{{\cal G}}

\def\su{\subset}

\def\se{\setminus}

\def\al{\alpha}
\def\be{\beta}

\def\De{\Delta}

\def\si{\sigma}

\def\La{\Lambda}

\def\cd{\cdot}
\def\stb{,\ldots ,}

\def\msk{\medskip}
\def\bsk{\bigskip}
\def\noi{\noindent}

\def\sumin{\sum_{i=1}^n}

\def\sumjn{\sum_{j=1}^n}
\def\sumjm{\sum_{j=1}^m}
\def\sumj0m{\sum_{j=0}^m}

\def\sumjk{\sum_{j=1}^k}
\def\sumi0n{\sum_{i=0}^n}

\begin{tiny}\end{tiny}

\def\proof{\noi {\bf Proof.} }

\def\x1n{x_1 \stb x_n}
\def\y1n{y_1 \stb y_n}

\def\msk{\medskip}
\def\bsk{\bigskip}
\def\noi{\noindent}

\def\stb{,\dots,}
\def\sik{\mathbb{R}^2}

\def\sumin{\sum_{i=1}^n}

\def\sumjn{\sum_{j=1}^n}
\def\sumjm{\sum_{j=1}^m}
\def\sumj0m{\sum_{j=0}^m}

\def\sumjk{\sum_{j=1}^k}

\newcommand{\R}{\mathbb{R}}
\newcommand{\C}{\mathbb{C}}
\newcommand{\Z}{\mathbb{Z}}

\begin{document}

\title{The discrete Pompeiu problem on the plane}

\author{Gergely Kiss, Mikl\'os Laczkovich and Csaba Vincze}

\footnotetext[3]{G. Kiss is partially supported by the project
R-AGR-0500-MRO3 of the University of Luxembourg, and also by the
Hungarian National Research, Development and Innovation Office,
Grant No. NKFIH 104178.} \footnotetext[4]{M. Laczkovich is partially
supported by the Hungarian National Research, Development and
Innovation Office, Grant No. NKFIH 104178.} \footnotetext[5]{Cs.
Vincze is supported by the University of Debrecen's internal
research project RH/885/2013.}

\maketitle

\begin{abstract}
We say that a finite subset $E$ of the Euclidean plane $\R^2$ has
the discrete Pompeiu property with respect to isometries (similarities),
if, whenever $f:\R^2\to \C$ is such that the sum of the values of $f$
on any congruent (similar) copy of $E$ is zero, then $f$ is identically
zero. We show that every parallelogram and every quadrangle with
rational coordinates has the discrete Pompeiu property w.r.t.
isometries. We also present a family of quadrangles depending on a continuous
parameter having the same property. We investigate the weighted version
of the discrete Pompeiu property as well, and show that every finite linear set
with commensurable distances has the weighted discrete Pompeiu property w.r.t.
isometries, and every finite set
has the weighted discrete Pompeiu property w.r.t. similarities.
\end{abstract}

\section{Introduction}\label{s1}
Let $K$ be a compact subset of the plane having positive
Lebesgue measure. The set $K$ is said to have the
Pompeiu property if the following condition is satisfied:
whenever $f$ is a continuous function defined on the plane, and
the integral of $f$ over every congruent copy of $K$ is zero, then
$f\equiv 0$. It is known that the closed disc does not have
the Pompeiu property, while all polygons have.
(As for the history of the problem, see \cite{R} and \cite{Z}.)

Replacing the Lebesgue measure by the counting measure, and the
isometry group by an arbitrary family $\gg$ of bijections mapping a set $X$
onto itself, we obtain the following notion. Let $K$ be a nonempty finite
subset of $X$. We say that {\it $K$ has the discrete Pompeiu property with
respect to the family $\gg$
if the following condition is satisfied:
whenever $f\colon X\to \cc$ is such that $\sum_{x\in K} f(\phi (x))=0$
for every $\phi \in \gg$, then $f\equiv 0$.}

We also introduce the weighted version of the discrete Pompeiu property.
We say that the $n$-tuple {\it $K = (x_1 \stb x_n )$
has the weighted discrete Pompeiu property with respect to the family $\gg$
if the following condition is satisfied: whenever $\al _1 \stb \al _n$ are
complex numbers with $\sumin \al _i \ne 0$ and $f\colon X\to \cc$ is such
that $\sumjn a_j f(\phi (x_j ))=0$ for every $\phi \in \gg$, then $f\equiv 0$.}

Apparently, the first results concerning the discrete Pompeiu property
appeared in \cite{Ze}, where the author considers
the Pompeiu problem for finite subsets of $\zz ^n$ w.r.t. translations.
The interest in the topic revived shortly
after the 70th William Lowell Putnam Mathematical
Competition (2009), where the following problem was posed:
Let $f$ be a real-valued function on the plane such that for every
square $ABCD$ in the plane, $f(A)+f(B)+f(C)+f(D)=0.$ Does it
follow that $f\equiv 0$? This is nothing but asking whether the set of vertices
of a square has the discrete Pompeiu property with respect to the similarities
of the plane. This problem motivated the paper \cite{GD}
by C. de Groote and M. Duerinckx. They prove that
every finite and nonempty subset of $\sik$ has
the discrete Pompeiu property w.r.t. direct similarities.
Another generalization of the Putnam problem appeared in \cite{KKS},
where it is proved that the set of vertices
of a square has the discrete Pompeiu property with respect to the group of
isometries. Recently, M. J. Puls \cite{Pu}
considered the discrete Pompeiu problem in groups.

In this paper we improve the results of \cite{GD} and \cite{KKS}.
We show that every finite and nonempty subset of $\sik$ has
the weighted discrete Pompeiu property w.r.t. direct similarities
(Theorem \ref{t1}).
We also show that the set of vertices of every parallelogram
has the discrete Pompeiu property with respect to the group of rigid motions
(Theorem \ref{t3}). We show the same for quadrangles with
rational coordinates (Theorem
\ref{tfourrac}), and for a family of quadrangles depending on a
continuous parameter (Theorem \ref{t4}). We also prove that in $\sik$ all
linear sets with commensurable distances have
the discrete Pompeiu property w.r.t. rigid motions
(Theorem \ref{t2}).

These results motivate the following questions: is it true that every
four element subset of the plane has the
discrete Pompeiu property with respect to the group of isometries?
Is it true that every nonempty and finite subset of the plane has the
same property? We do not know the answer.

We conclude this introduction with a remark concerning the
family of translations in an Abelian group. As the following proposition
shows, this family is `too small': finite sets, in general, cannot have
the discrete Pompeiu property w.r.t. this group.

\begin{proposition}\label{p1}
Let $G$ be a torsion free Abelian group. If $E$ is a finite subset of $G$
containing at least 2 elements, then $E$ does not have the discrete
Pompeiu property w.r.t. the family of all translations of $G$.
\end{proposition}

\proof
Note that if the torsion free rank of $G$ is less than continuum, then
this is a special case of \cite[Theorem 3.1]{Pu}. In the general case
let $H$ be the subgroup of $G$ generated by $E$. Then
$H$ is a finitely generated torsion free Abelian group, and thus $H$
is isomorphic to $\Z^n$ for some finite $n$.
By Zeilberger's theorem \cite{Ze}, $E$ does not have the discrete
Pompeiu property in $H$ w.r.t. the family of translations; that is,
there is a nonzero function $f\colon H\to \cc$ such that the sum of
the values of $f$ taken on any translate of $E$ is zero.
It is clear that we can find such a function on every coset of $H$.
The union of these functions has the same property on $G$, showing that
$E$ does not have the discrete
Pompeiu property in $G$ w.r.t. the family of translations on $G$.
$\square$

\msk
In the proposition above we cannot omit the requirement that
$G$ be torsion free. E.g., if $G$ is a finite group having $n\ge 3$
elements and $E$ is a subset of $G$ having $n-1$ elements,
then $E$ has the discrete Pompeiu property w.r.t. translations. Indeed,
if the sum of the values of $f$ is zero
on each translate of $E$ then $f$ must be constant, and the constant must
be zero.

\section{Preliminaries: generalized polynomials and
exponential functions on Abelian groups}\label{s2}

Let $G$ be an Abelian group. If $f\colon G\to \cc$ and $h\in G$,
then $\De _h f$ denotes the function defined by $\De _h f (x)=f(x+h )-f(x)$
$(x\in G)$. The function $f\colon G\to \cc$ is said to be a
{\it generalized polynomial} if there is an $n$ such that
$\De _{h_1} \ldots \De _{h_{n+1}} f \equiv 0$ for every ${h_1} \stb {h_{n+1}} \in
G$. The degree of $f$ is the smallest such $n$. Thus the
generalized polynomials of degree zero are the nonzero constant functions.
The degree of the identically zero function is $-1$ by definition.

The function $g\colon G\to \cc$ is an {\it exponential}, if $g\ne 0$ and $g(x+y)=g(x)\cd g(y)$ for every $x,y\in G$.
By a {\it monom} we mean a function of the form $p\cd g$, where $p$ is a
generalized polynomial, and $g$ is an exponential. Finite sums of monoms are
called {\it polynomial-exponential functions}.

Let $\C^G$ denote the linear space of all complex valued functions defined
on $G$ equipped with the product topology.
By a {\it variety on $G$} we mean a translation invariant, closed, linear
subspace of $\C^G$.
We say that spectral analysis holds in $G$, if every nonzero
variety contains an exponential function.

We shall need the fact that spectral analysis holds in every
finitely generated and torsion free Abelian group. In fact,
this is true in every Abelian group whose torsion free rank is less
than continuum \cite{LS}. However, for finitely generated and torsion
free Abelian groups this also follows from Lefranc's theorem.
Lefranc proved in \cite{LF} that if $n$ is finite then
{\it spectral synthesis} holds in $\zz ^n$; that is, every variety on $\zz ^n$ is spanned
by polynomial-exponential functions. Therefore, if a variety $V$ on $\zz ^n$
contains nonzero functions, then it has to contain nonzero
polynomial-exponential functions. It is easy to see that if a
polynomial-exponential function $\sumin p_i \cd g_i$ is contained
in a variety $V$, where $p_1 \stb p_n$ are nonzero generalized polynomials
and $g_1 \stb g_n$ are distinct exponentials, then necessarily $g_i \in V$
holds for every $i=1\stb n$. Since every
finitely generated and torsion free Abelian group
is isomorphic to $\zz ^n$ for some finite $n$, it follows that
spectral analysis (and, in fact, spectral synthesis) holds in such groups.
We shall need the following special case.

\begin{lemma}\label{l2}
Let $G$ be a finitely generated subgroup of the additive group of $\cc$,
let $\al _{j,k}, b_{j,k}$ $(j=1\stb n, \ k=1\stb m)$ be complex numbers,
and let $Y$ be a subset of $\cc$. Let $V$ denote the set of functions
$f\colon G\to \cc$ such that
$$\sumjn \al _{j,k} f(x+b_{j,k} y)=0$$
for every $k=1\stb m$, $x\in G$ and $y\in Y$ satisfying $b_{j,k} y \in G$
for every $j=1\stb n$ and $k=1\stb m$. Then $V$ is a variety on $G$.
Consequently, if $V$ contains a non-identically zero function, then
$V$ contains an exponential function defined on $G$.
\end{lemma}

\proof
It is clear that $V$ is a translation invariant linear
subspace of $\C^G$. Since $G$ is countable, the topology of
$\C^G$ is the topology determined by pointwise convergence.
Obviously, if $f_i \in V$ and $f_i \to f$ pointwise on $G$, then
$f\in V$. Thus $V$ is closed. $\square$

\section{Similarities}\label{s3}

It was shown by C. de Groote and M. Duerinckx in
\cite{GD} that {\it every finite and nonempty subset of $\sik$ has
the discrete Pompeiu property w.r.t. direct similarities.} By a
direct similarity we mean a transformation that is a composition
of translations, rotations and homothetic transformations. The authors
also discuss the possible generalizations when $\sik$ is replaced
by $K^p$ where $K$ is a field, and the transformation group is a
subgroup of $AGL(p,K)$. We note that the argument given by
C. de Groote and M. Duerinckx also proves the following generalization.

\begin{proposition}\label{p2}
Let $\gg$ be a transitive and locally commutative
transformation group acting on $X$
such that for every $x,y,z\in X$ with $y\ne x\ne z$ there exists a map
$f\in \gg$ such that $f(x)=x$ and $f(y)=z$. Then every finite and nonempty
proper subset of $X$ has the discrete Pompeiu property w.r.t. $\gg$.
\end{proposition}

We say that a transformation $g\colon \rr \to \rr$ is an
order preserving similarity, if $g(x)=a+cx$ for every $x\in \rr$, where
$a\in \rr$ and $c>0$.
\begin{proposition}\label{p3}
Every finite and nonempty subset
of $\rr$ has the discrete Pompeiu property w.r.t. the group of order
preserving similarities.
\end{proposition}
\proof
Although Proposition \ref{p2} cannot be applied directly, a variant of
the argument given by C. de Groote and M. Duerinckx in \cite{GD} works.
Let $E=\{x_1 \stb x_n \}$. Suppose that $f\colon \rr \to \rr$ is such that
$\sumin f(a+cx_i )=0$ for every $a\in \rr$ and $c>0$.
Replacing $E$ by a translated copy we may assume
that $0=x_1 <x_2 < \ldots <x_n$. We put $A_i =\{ x_i +x_i x_j \colon j=2\stb n\}$
and $B_j =\{ x_i +x_i x_j \colon i=2\stb n\}$. Then $A_i \cup \{ x_i \}$
is the image of $E$ under an order preserving similarity for every $i\ge 2$,
and thus $\sum_{j=2}^n f(x_i +x_i x_j )=-f(x_i )$ $(i=2\stb n)$. Similarly,
$B_j \cup \{ 0 \}$ is the image of $E$ under an order preserving
similarity for every $j\ge 2$, and thus
$\sum_{i=2}^n f(x_i +x_i x_j )=-f(0 )$ $(j=2\stb n)$. Therefore,
\begin{align*}
f(0)= & -\sum_{i=2}^n f(x_i )= \sum_{i=2}^n \sum_{j=2}^n f(x_i +x_i x_j ) =
\sum_{j=2}^n \sum_{i=2}^n f(x_i +x_i x_j ) = \\
&=\sum_{j=2}^n (-f(0))=-(n-1)f(0).
\end{align*}
Thus we have $f(0)=0$. For every $b\in \rr$, the function $T_b f$
defined by $T_b f (x)=f(x+b)$ also satisfies the condition
$\sumin T_b f(a+cx_i )=0$ for every $a\in \rr$ and $c>0$. Therefore,
$T_b f(0)=f(b)=0$ for every $b\in \rr$. $\square$

\msk
De Groote and M. Duerinckx ask in \cite{GD} if the finite subsets of
the plane have the weighted discrete Pompeiu property w.r.t.
direct similarities. In the next theorem we show that the answer is affirmative.

\begin{theorem}\label{t1}
Every $n$-tuple of distinct points of $\sik$ has the weighted discrete
Pompeiu property w.r.t. direct similarities.
\end{theorem}

\proof
We identify $\sik$ with the complex plane $\cc$. We put
$\cc ^* =\cc \se \{ 0\}$.

Let $(b_1 \stb b_n )$ be an $n$-tuple of distinct complex numbers.
Let $\al _1 \stb \al _n$ be complex numbers such that $\sumin \al _i \ne 0$,
and let $f\colon \cc \to \cc$ be such that
\begin{equation}\label{e2}
\sumin \al _i f(x+b_i y)=0
\end{equation}
for every $x\in \cc$ and $y\in \cc ^*$.
We have to prove that $f\equiv 0$.

If \eqref{e2} holds for every $x,y\in \cc$, then $f\equiv 0$ is one of
the statements of \cite[Theorem 2.4]{KV}. Therefore, it is
enough to show that if
\eqref{e2} holds for every $x\in \cc$ and $y\in \cc ^*$, then
it holds for every $x,y\in \cc$. In the following theorem
we shall prove more.

We say that a family $I$ of subsets of $\cc$ is a proper
and translation invariant ideal, if
$A,B\in I$ implies $A\cup B\in I$, $A\in I$ and $B\su A$ implies $B\in I$,
$\cc \notin I$, and if $A\in I$ then $A+c =\{ x+c \colon x\in A\} \in I$
for every $c\in \cc$. It is clear that the family of finite subsets
of $\rr$ is a proper and translation invariant ideal.

\begin{theorem}\label{t5}
Let $I$ be a proper and translation invariant
ideal of subsets of $\C$.
Let $b_1 \stb b_n$ be distinct complex numbers, and suppose that
the functions $f_1 \stb f_n \colon \C \to \C$ satisfy
\begin{equation}\label{el1}
\sumin f_i (x+b_i y)=0
\end{equation}
for every $x\in \C$ and $y\in \cc \se A$, where $A\in I$.
Then each $f_i$ is a generalized polynomial of degree $\le n-2$,
and \eqref{el1} holds for every $x,y\in \cc$.
\end{theorem}

\proof
First we prove that each $f_i$ is a generalized polynomial of degree $\le n-2$.
We prove by induction on $n$. The case of $n=1$ is obvious.

Now let $n\ge 2$, and suppose that the statement is true for $n-1$.
Let $f_1 \stb f_n$ satisfy \eqref{el1} for every
$x$ and $y\notin A$, where $A\in I$. Since the role of the functions $f_i$
is symmetric, it is enough to prove that $f_1$ is
a generalized polynomial of degree $\le n-2$.
Note that $b_1 \ne b_n$ by assumption. Let $h\in \C$
be fixed. Then we have
\begin{equation}\label{e3}
\sum_{i=1}^{n} f_i (x+h+b_i y) =0
\end{equation}
for every $x$ and $y\in \cc \se A$, and
\begin{equation}\label{e4}
\sum_{i=1}^{n} f_i (u+b_i v ) =0
\end{equation}
for every $u$ and $v\notin A$. Substituting $u=x-b_1 h/(b_n -b_1 )$
and $v=y+h/(b_n -b_1 )$ into \eqref{e4} and subtracting from \eqref{e3}
we obtain
$$\Delta _h f_1 (x+b_1 y) +\sum_{i=2}^{n-1} \left[ f_i (x+h+ b_i y) -
f_i \left( x+ \frac{b_i -b_1}{b_n -b_i} h +b_i y \right) \right] =0$$
for every $y$ such that $y\notin A$ and $v=y+h/(b_n -b_1 ) \notin A$.
(If $n=2$ then the sum on the left hand side is empty.)
Putting $g_i (z)=f_i (z+h) - f_i \left(
z+ \tfrac{b_i -b_1}{b_n -b_i}h \right)$ $(z\in \C )$,
we obtain that
$$\Delta _h f_1 (x+b_1 y) +\sum_{i=2}^{n-1} g_i (x+b_i y) =0$$
for every $x$ and for every $y\notin A\cup (A-h/(b_n -b_1 ))$.
Since $A\cup (A-h/(b_n -b_1 ))\in I$, it follows from the induction
hypothesis that $\Delta _h f_1$ is a
generalized polynomial of degree $\le n-3$.
As this is true for every $h$, we obtain that $f_1$
is a generalized polynomial of degree $\le n-2$.

We still have to prove that \eqref{el1} holds for every $x,y\in \cc$.
Let $x\in \cc$ be fixed, and put $G(y)=\sum_{i=1}^{n} f_i (x+b_i y)$
for every $y\in \cc$. We have to prove that $G(y)=0$ for every
$y\in \cc$.

It is easy to see that if $f$ is a generalized polynomial,
then so is $y\mapsto f (x+b y) =g(y)$. This can be proved
by induction on the degree of $f$, using $\De _h g (y)=\De _{bh} f(x+by)$.
Since each $f_i$ is a generalized polynomial, it follows
that so is $g_i (y)=f_i (x+b_i y)$ for every $i$, and thus so is
$G=g_1 +\ldots +g_n$.

We know that $G(y)=0$ for every $y\notin A$. Therefore, in order to
prove $G\equiv 0$, it is enough to show that if $f\colon \cc \to \cc$
is a generalized polynomial and $f(x)=0$ for every $x\in \cc \se A$ where
$A\in I$, then $f\equiv 0$. We prove by induction on the degree of $f$.
The statement is obvious if the degree is $\le 0$. Indeed, in this case
$f$ is constant, and has a value equal to zero, since $I$ is a proper ideal.
Suppose the degree of $f$ is $n>0$, and the statement
is true for generalized polynomials of degree $<n$. For every $h$,
we have $\De _h f (x)=0$ for every $x\in \cc \se (A \cup (A-h))$. Since
$A \cup (A-h) \in I$, it follows from the induction hypothesis that
$\De _h f(x)=0$ for every $x$. This is true for every $h$,
which shows that $f$ is constant. As we saw above, the constant must be zero.
This completes the proof.
\hfill $\square$

\section{Isometries and rigid motions: some general remarks}\label{s4}

By a rigid motion we mean an isometry that preserves orientation.
An isometry of $\sik$ is a rigid motion if it is a translation or a rotation.
\begin{proposition}\label{p5}
Every subset of the plane containing $1$, $2$ or $3$ points has
the discrete Pompeiu property w.r.t. rigid motions.
\end{proposition}

\proof
The case of the singletons is obvious.
Let $E=\{a,b\}$ and $r=|a-b| >0$. Suppose that $f\colon \sik \to \cc$
is such that
$f(\si (a))+f(\si (b))=0$ for every rigid motion $\si$. Then $f$ has the same
value at every pair of points $a_1,a_2$ with distance $\le 2r$.
Indeed, there is a point $b$ such that $|b-a_i| =r$ $(i=1,2)$,
and thus $f(a_1 )=-f(b)=f(a_2 )$. Now, any two points $a,b \in \sik$ can be
joined by a sequence of points $a=a_0 \stb a_n =b$ such that $|a_i -a_{i-1}| \le 2r$, and thus $f(a)=f(b)$. Therefore, $f$ must be constant, and the value of
the constant must be zero.

Let $H=\{a,b,c\}$, where $a,b,c$ are distinct, and let $f\colon \sik \to \cc$
be such that $f(\si (a))+f(\si (b)) +f(\si (c))=0$ for every rigid motion $\si$.
By changing the notation of the points $a,b,c$ if necessary,
we may assume that $c\ne (a+b)/2$. Let
$c'=a+b-c$. Then $c'$ is the
reflection of $c$ about middle point of the segment $[a,b]$, and thus
$f(\si (b))+f(\si (a)) +f(\si (c'))=0$ for every rigid motion $\si$.
Thus $f(\si (c'))=f(\si (c))$ for every rigid motion $\si$, which implies that
$f(x)=f(y)$ whenever $|x-y|=|c' -c|$. The argument above shows that
$f$ is constant, and, in fact, $f\equiv 0$. $\square$

\begin{remark} \label{r1}
{\rm It is easy to see that if $n\le 2$,
then every $n$-tuple has the weighted discrete Pompeiu property w.r.t.
isometries. The same is true for those triplets
$(a,b,c)$ whose points are not collinear. In this case we have to
modify the proof above by choosing the point $c'$ to be the reflection of
$c$ about the line going through $a$ and $b$ instead of the point $a+b-c$
in order to avoid changing the weights of $a$ and $b$.}
\end{remark}

\begin{proposition}\label{p6}
Let $E$ be a finite set in the plane.
If there exists an isometry $\sigma$ such that $|E\cap\sigma(E)|=|E|-1$, then
$E$ has the discrete Pompeiu property w.r.t. isometries.
\end{proposition}

\proof Let $E\se \si (E)=\{ a\}$ and $\si (E) \se E=\{ b \}$.
If $f\colon X\to \cc$ is such that $\sum_{x\in \phi  (E)}  f(x)=0$
for every isometry $\phi$, then taking the difference
of the equations $\sum_{x\in (\phi \si ) (E)}  f(x)=0$ and
$\sum_{x\in \phi  (E)}  f(x)=0$, we obtain
$f(\phi (a))=f(\phi (b))$ for every isometry $\phi$.
Thus $f(x)=f(y)$ whenever $|x-y|=|a-b|$. As we saw before, this implies that
$f$ is identically zero. $\square$

\begin{remark} Concerning the discrete Pompeiu property in higher dimensions,
we note that Proposition \ref{p6} holds without any essential modification
in $\mathbb{R}^n$ for every $n\ge 2$.
As for Proposition \ref{p5}, it is easy to see that
every subset of $\mathbb{R}^n$ $(n\ge 2)$ containing affinely independent
points has the discrete Pompeiu property w. r. t. isometries. Using an
inductive argument it is enough to consider the case of $n+1$ points
in general position. Such a set satisfies the condition of
Proposition \ref{p6}: let $\si$ be the reflection about a facet.
\end{remark}

\bsk
By Proposition \ref{p6}, if a set $E$ consists of consecutive
vertices of a regular $n$-gon $R$, and $E\ne R$, then $E$
has the discrete Pompeiu property w.r.t. isometries. Also, if $E$ is a finite
set of collinear points forming an arithmetic progression, then
$E$ has the discrete Pompeiu property w.r.t. isometries.
Our following theorem is the generalization of this fact.

\begin{theorem} \label{t2}
Let $E$ be an $n$-tuple of collinear points in $\sik$
with pairwise commensurable distances. Then
$E$ has the weighted discrete Pompeiu property w.r.t. rigid motions
of the plane.
\end{theorem}

\begin{lemma} \label{l5}
Let $x_1 \stb x_n ,y_1 \stb y_k \in \sik$
and $\al _1 \stb \al _n ,\be _1 \stb \be _k \in \cc$ be such that
\begin{itemize}
\item[{\rm (i)}]
$y_1 \stb y_k$ are collinear with commensurable distances,
\item[{\rm (ii)}]
$\sumin \al _i \ne 0$, and
\item[{\rm (iii)}]
at least one of the numbers $\be _1 \stb \be _k$ is nonzero.
\end{itemize}
If $f\colon \sik \to \rr$ is such that
\begin{equation}\label{e24}
\sumin \al _i f(\si (x_i ))=\sumjk \be _j f(\si (y_j ))=0
\end{equation}
for every rigid motion $\si$, then $f$ is identically zero.
\end{lemma}

\proof
We identify $\sik$ with the complex plane $\cc$. We put
$\cc ^* =\cc \se \{ 0\}$ and $S^1 =\{ x\in \cc \colon |x|=1\}$.
Then every rigid motion is of the
form $x\mapsto a+ux$ $(x\in \cc )$, where $a\in \cc$ and $u\in S^1$.

Let $a,c\in \cc$ and $c\ne 0$. If we replace $x_i$ by $a+cx_i$,
$y_j$ by $a+cy_j$ for every $i$ and $j$, and replace
$f$ by $f_1 (x)=(x/c)$, then \eqref{e24} remains valid
for every rigid motion $\si$. Indeed, for every $\si$, the map
$x\mapsto \si (a+cx) /c$ is a rigid motion if and only if $\si$ is.
Note that if $f_1$ is identically zero, then so is $f$.

Therefore, replacing $x_i$ by $a+cx_i$,
$y_j$ by $a+cy_j$ for every $i=1\stb n$ and $j=1\stb k$ with a suitable
$a\in \cc$ and $c\in \cc ^*$, we may assume that $y_1 \stb y_k$
are positive integers.
By supplementing the system if necessary, we may assume that
$y_j =j$ $(j=1\stb m)$. We put $\be _j =0$ for every added $j$. Then we have
\begin{equation}\label{e1}
\sumin \al _i f (x+ ux_i)= \sumjm \be _j f(x+ju )=0
\end{equation}
for every $x\in \C$ and $u\in S^1$. We show that this implies $f\equiv 0$.
Suppose that $f$ is not identically zero, and let $z_0 \in \C$ be such that
$f(z_0 )\ne 0$.

Let $K$ be an integer greater than $\max _{1\le i\le n} |x_i |$.
It is clear that every $z\in \cc$ with $|z|<K$ is the
sum of $K$ elements of $S^1$. Let
$U$ be a finite subset of $S^1$ such that $1\in U$, and
$x_i /\nu$ is the sum of
$K$ elements of $U$ for every $i=1\stb n$ and
$\nu =1\stb N$, where $N=m^{K\cd m^K}$.

Let $G$ denote the additive subgroup of $\C$ generated by the
elements $z_0$, $u\in U$ and $u x_i$ $(u\in U, \ i=1\stb n)$.
Then $G$ is a finitely generated subgroup of $\C$.
Let $V$ denote the set of functions defined on $G$ and
satisfying \eqref{e1} for every $x\in G$ and $u\in U$.
The set of functions $V$ contains the restriction of $f$ to $G$ which
is not identically zero, as $z_0 \in G$.
Therefore, by Lemma \ref{l2}, $V$ contains an
exponential function $g \colon G\to \C$.

If $u\in U$, then \eqref{e1} gives $\sumjm \be _j g (u)^j =0$.
Therefore, $g (u)$ is a root of the polynomial
$p(x)=\sum_{j=1}^{m} \be _j x^{j-1}$.
Let $\Lambda$ denote the set of the nonzero roots of $p$.
Then $\Lambda$ has at most $m-1$ elements,
and $g (u)\in \Lambda$ for every $u\in U$. For every
$i=1\stb n$ and $\nu =1\stb N$, $x_i /\nu$ is the sum of $K$
elements of $U$. Thus $g(x_i /\nu )$ is the product of $K$ elements
of $g(U)\su \La$. Therefore, the set $F=\{ g(x_i /\nu ) \colon i=1\stb n, \
\nu =1\stb N\}$ has less than $m^K$ elements.

Let $1\le i\le n$ be fixed, and put $g (x_i)=c$. We prove $c=1$. Since
$g(x_i /\nu )\in F$ for every $\nu =1\stb m^K$, there are integers
$1\le \nu <\mu \le m^K$ such that $g (x_i /\nu  )=g (x_i /\mu )$.
Then
$$c^\mu  =g (x_i /\nu )^{\nu \mu}=g (x_i /\mu )^{\nu \mu} =c^\nu ,$$
and thus $c^{\mu -\nu }=1$. Let $\mu -\nu =s$, then  $s<m^K$ and $c^s =1$.
If $s=1$, then $c=1$ is proved. If $s>1$, then, by $g (x_i /s^t )\in F$
for every $t=1\stb m^K$, there are integers
$1\le r<t\le m^K$ and there is an element $b\in F$ such that
$g (x_i /s^r )=g (x_i /s^t )=b$. Then
$$c=g(x_i )=b^{s^t} =b^{s^r \cd s^{t -r}} =c^{s^{t-r}} =1,$$
since $c^s =1$. This proves $g (x_i)=1$ for every $i=1\stb n$.

Then, applying \eqref{e1} with $x=0$ and $u=1$, we obtain $\sumin \al _i =0$
which is impossible. This contradiction completes the proof.
\hfill $\square$

\msk \noi
{\bf Proof of Theorem \ref{t2}.} Let $E=(x_1 \stb x_n )$, where
$x_1 \stb x_n$ are collinear with commensurable distances. Let $\al _1
\stb \al _n$ be complex numbers with $\sumin \al _i \ne 0$, and let $f\colon \cc \to \cc$ satisfy $\sumin \al _i f(\si (x_i ))=0$ for every rigid motion $\si$.
Applying Lemma \ref{l5} with $k=n$, $y_i =x_i$ and $\be _i =\al _i$
$(i=1\stb n)$, we obtain that $f$ is identically zero. \hfill $\square$

\begin{remark}
{\rm The isometry group of $\rr$ consists of translations and reflections.
Since no finite subset of $\rr$ has the discrete Pompeiu property
w.r.t. translations by Proposition \ref{p1}, and every reflected copy of the set $\{ 1\stb n\}$
is also a translated copy, it follows that} the set $\{ 1\stb n\}$
does not have the discrete Pompeiu property
w.r.t. isometries of $\rr$. {\rm (This is why we had to
step out from $\rr$ into the plane in the proof of Theorem \ref{t2}.)
Note, however, that there are subsets of $\zz$ which have
the discrete Pompeiu property w.r.t. isometries of $\rr$.
The set of integers $0=z_0 < z_1 < \ldots < z_k$
has this property if and only if the polynomials
$p(x)=\sum_{i=0}^k x^{z_i}$ and $q(x)=\sum_{i=0}^k x^{z_k - z_i}$ have no common
roots. (This follows immediately from Zeilberger's theorem \cite{Ze}.)
This condition is clearly satisfied if the set
$\{z_0, z_1, \ldots, z_k\}$ is not symmetric about the point $(z_0+z_k )/2$,
and if $p$ is irreducible in $\mathbb{Z}[x]$.

Since each coefficient of $p$ is $0$ or $1$, it is easy to decide
whether $p$ is irreducible or not. If there is an $n\ge 3$ such that
$p(n)$ is prime, then $p$ is irreducible (see \cite{M}).
By the Buniakowski-Schinzel conjecture, this condition is also necessary
for the irreducibility of $p$.
}
\end{remark}

\section{Quadrangles under isometries}\label{s5}

\begin{theorem} \label{t3}
The set of vertices of any parallelogram has the discrete Pompeiu property
w.r.t. rigid motions of the plane.
\end{theorem}

\proof
We identify $\sik$ with the complex plane $\cc$. We put
$\cc ^* =\cc \se \{ 0\}$ and $S^1 =\{ u\in \cc \colon |u|=1\}$.
Let $E$ be a set of vertices of a parallelogram.
Without loss of generality we may assume that $0\in E$.
Then $E=\{ 0,a,b,a+b\}$, where $0 \ne a,b\in \C$ and $a\ne b$. Clearly, it
is enough to prove that if $f\colon \cc \to \cc$ is such that
\begin{equation}\label{epo4}
f(x)+f(x+ay)+f(x+by)+f(x+(a+b)y)=0
\end{equation}
for every $x\in \C$ and $y\in S^1$, then $f\equiv 0$.
Suppose that there exists a nonzero $f$ satisfying \eqref{epo4}, and
let $z_0\in \C$ be such that $f(z_0)\ne 0$.

Let $F$ be a finite subset of $\cc$, and let $G$ denote the additive
subgroup of $\cc$ generated by $F\cup \{ z_0\}$.
Let $V$ denote the set of functions $f\colon G\to \C$ satisfying
\eqref{epo4} for every $x\in G$ and $y\in S^1_G =\{ y\in S^1 \colon
ay,by\in G \}$. Since $f|_G \in V$ and $z_0 \in G$, it follows that $V\ne 0$.

By Lemma \ref{l2}, there exists an exponential
function $g$ in $V$. Since $g$ satisfies \eqref{epo4} and
$g(x+ay)=g(x)g(ay)$ and $g(x+(a+b)y) =g(x)g(ay) g(by)$, we obtain
$g(x)(1+g(ay)+g(by)+g(ay)g(by))=0$
whenever $x\in G$ and $y\in S^1_G $.
Since $g(x)\ne 0$, we get $(1+g(ay))(1+g(by))=0$
for every $y\in  S^1 _G$ . That is, we have either
$g(ay)=-1$ or $g(by)=-1$ for every $y\in S^1_G$.

Let $P$ be an arbitrary parallelogram obtained from $E$ by a rigid motion
and having vertices in $G$. Then the vertices of $P$ are $c=x$,
$d=x+ay$, $e=x+(a+b)y$, $f=x+by$ with a suitable $x\in G$ and $y\in S^1 _G$.
Then we have either
$g(d)/g(c)=g(e)/g(f)=-1$ or $g(f)/g(c)=g(e)/g(d)=-1$.
In other words, the values of $g$ at the points $c,d,e,f$ are either
$g(c), -g(c),g(e), -g(e)$ or $g(c), g(d), -g(d), -g(c)$.
Therefore, the vertex set of the parallelogram can be decomposed into two pairs
with $g$-values of the form $(x,-x)$ in each pair.

Let ${\cc}^*= X_1 \cup X_2$ be a decomposition of ${\cc}^*$ such that
$X_1=-X_2$. Let $h(x)=1$ if $g(x)\in X_1$, and $h(x)=-1$ if $g(x)\in X_2$.
Then $h\colon G \to \{1,-1\}$ has the following property: if $\si$
is a rigid motion and if $\si (E)\su G$, then there are two elements of
$\si (E)$
where the function $h$ takes the value $1$, and at the other two elements
of $\si (E)$ the function $h$ takes the value $-1$.

Since this is true for every group generated by any finite subsets of $\sik$,
we may apply Rado's selection principle \cite{G}. We find that there
exists a function $h\colon \sik \to \{ 1,-1\}$ such that whenever $\si$
is a rigid motion, then there are two elements of $\si (E)$
where the function $h$ takes the value $1$, and at the other two elements
of $\si (E)$ the function $h$ takes the value $-1$.

The existence of such a function, however, contradicts a known fact
of Euclidean Ramsey theory.
By a theorem of Shader \cite[Theorem 3]{Sh}, for every $2$-coloring of the
plane, and for every parallelogram $E$, there is a congruent copy $P$
of $E$ such that at least three vertices of $P$ has the same color.
It is clear from the proof that $P$ can be obtained from $E$ by a rigid motion.
(See the Remark on p. 563 in \cite{EG}.)
This contradicts the existence of the function $h$ with the properties
described, proving that $f$ must be identically zero. $\square$

\bsk
Our next aim is to prove the following.
\begin{theorem}\label{tfourrac}
Every set $E\su \sik$ of four points having rational
coordinates has the weighted discrete Pompeiu property w.r.t. the
group of isometries of $\sik$.
\end{theorem}

\proof
If the points of $E$ are collinear, then the statement is a consequence of
Theorem \ref{t2}. If there are three collinear points of $E$,
then the statement follows from Proposition \ref{p6}. Therefore, we may assume
that the points of $E$ are in general position.
Let $E=(x_1 \stb x_4 )$. By changing the order of the indices
we may assume that $x_1$ and $x_2$ are vertices of the convex hull of $E$.

Let $\al _1 \stb \al _4$ be complex numbers such that
$\sum _{j=1}^4 \al _j \ne 0$, and let $f\colon \cc \to \cc$ be such that
\begin{equation} \label{e18}
\sum_{j=1}^4 \al _j f(\si (x_j ))=0
\end{equation}
for every isometry $\si$. We have
to show that $f$ is identically zero.
If any of the numbers $\al _1 \stb \al _4$ is zero then
$f\equiv 0$ follows from Remark \ref{r1}. Therefore, we may assume that
$\al _4 \ne 0$.

Let $\si _1$ be the
reflection about the line $\ell _1$ going through the points $x_1 ,x_2$.
Let $y_1 =\si _1 (x_4 )$ and $x_5 =\si _1 (x_3 )$,
then $y_1$ and $x_5$ have rational coordinates. We have, for every $\si$,
$(\si \circ \si _1 ) (x_i )=\si (x_i )$ for $i=1,2$ ,
$(\si \circ \si _1 ) (x_3 )=\si (x_5 )$ and
$(\si \circ \si _1 ) (x_4 )=\si (y_1 )$. Therefore
$$
\al _1 f(\si (x_1))+\al _2 f(\si (x_1))+\al _3 f(\si (x_5 ))+\al _4 f(\si (y_1 ))
=0
$$
for every isometry $\si$. Subtracting \eqref{e18} we obtain
\begin{equation} \label{e7}
\al _3 f(\si (x_5 ))+\al _4 f(\si (y_1 )) -\al _3 f(\si (x_3 ))-
\al _4 f(\si (x_4 )) =0
\end{equation}
for every isometry $\si$. Suppose that the line going through the points
$x_3$ and $x_4$ is perpendicular to $\ell _1$. Then the points
$x_5 ,y_1 , x_3 ,x_4$ are collinear. They have rational coordinates, so
the distances between them are commensurable. Now $f$ satisfies
both \eqref{e18} and \eqref{e7} for every isometry $\si$, and thus,
by Lemma \ref{l5}, $f\equiv 0$.

Therefore, we may assume that the line going through the points
$x_3$ and $x_4$ is not perpendicular to $\ell _1$. Let $\si _2$ be the
reflection about the line $\ell _2$ going through the points $x_3 ,x_5$.
Note that the lines $\ell _1$ and $\ell _2$ are perpendicular.
We put $y_2 =\si _2 (y_1 )$,
$y_3 =\si _2 (x_4 )$ and $y_4 =x_4$. Then $y_1 ,y_2 ,y_3 ,y_4$
are the vertices of a rectangle $R$ listed either
clockwise or counter-clockwise.
It is clear that $y_1 ,y_2 ,y_3 ,y_4$ have rational coordinates.
We claim that
\begin{equation}\label{e14}
f(\si (y_1 )) -f(\si (y_2 )) +f(\si (y_3 )) -f(\si (y_4 )) =0
\end{equation}
holds for every isometry $\si$. Indeed,
$(\si \circ \si _2 ) (x_5 )= \si (x_5 )$,
$(\si \circ \si _2 )(x_3 )= \si (x_3 )$,
$(\si \circ \si _2 )(y_1 )= \si (y_2 )$ and
$(\si \circ \si _2 ) (x_4 )= \si (y_3 )$ and thus, by \eqref{e7} we obtain
$$\al _3 f(\si (x_5 ))+\al _4 f(\si (y_2 )) -\al _3 f(\si (x_3 ))-\al _4
f(\si (y_3 )) =0.$$
Subtracting \eqref{e7} and dividing by $-\al _4$ we obtain \eqref{e14}
for every isometry $\si$.

Since the coordinates of $y_1 \stb y_4$ are
rational, it follows that the side lengths of $R$ are commensurable.
(The side lengths themselves can be irrational.) Thus, there exists a square
$Q$ with vertices $z_1 \stb z_4$ such that $Q$ can be decomposed into
finitely many translated copy of $R$. If we add the equations \eqref{e14}
for those translations $\si$ that map $R$ into these translated copies,
then we get
\begin{equation}\label{e6}
f(z_1 ) -f(z_2 ) +f(z_3 ) -f(z_4 ) =0,
\end{equation}
since all other terms cancel out. By rescaling the set $E$ and also the function
$f$ if necessary, we may assume that the side length of $Q$ is $1$.
Clearly, \eqref{e6} must hold whenever $z_1 \stb z_4$ are the vertices
of a square of unit side length. That is, we have
\begin{equation*}
f(x) -f(x+u ) -f(x+u\cd i ) +f(x+u+u\cd i ) =0
\end{equation*}
for every $x\in \cc$ and $u\in S^1$.

Now we turn to the proof of $f\equiv 0$. Suppose this is not true, and
fix a $z_0 \in \cc$ such that $f(z_0 )\ne 0$. Let $a_1 \stb a_N$ be
vectors of length $12$ such that each of the numbers $x_1 \stb x_4$ is
the sum of some of the $a_j$'s. Let $u_j =a_j /12$ and
$v_j =(3u_j +4 u_j \cd i )/5$ for every $j=1\stb N$. Then $u_j, v_j$
are unit vectors for every $j$. Let $U$ denote the set of vectors
$$u_j, \, u_j \cd i, \,   v_j, \,  v_j \cd i  \qquad  (j=1\stb N),$$
and let $G$ denote the additive group generated by the set
$U\cup \{ x_j u\colon j=1\stb4,\ u\in U\} \cup \{ z_0 \}$.
Then $G$ is a finitely generated group.
Let $V$ be the set of functions $g\colon G\to \cc$
satisfying the following condition:
$$\sum_{j=1}^4 \al _j g(x+x_j \cd u)=0$$
and
\begin{equation}\label{e8}
g(x)-g(x+u)-g(x+u\cd i)+g(x+u+u\cd i )=0
\end{equation}
for every $x\in G$ and $u\in U$.

The set $V$ contains a non-identically zero function (namely the restriction of
$f$ to $G$), so by Lemma \ref{l2}, $V$ contains an exponential function $g$.
Then \eqref{e8} implies $(1-g(u))\cd (1-g(u\cd i))=0$, and thus
we have either $g(u)=1$ or $g(u\cd i)=1$ for every $u\in U$.

Now we show that $g(a_j )=1$ for every $j=1\stb N$.
If $g(u_j )=1$, then this follows from $g(a_j )=g(12 u_j )=g(u_j )^{12}$.
Therefore we may assume that $g(u_j \cd i)=1$.
Since $5v_j =3u_j +4u_j \cd i$ and $5v_j \cd i =-4u_j +3\cd u_j \cd i$,
we have $g(v_j )^5 =g(u_j )^3 \cd g(u_j \cd i)^4=g(u_j )^3$ and
$g(v_j \cd i )^5 =g(u_j )^{-4} \cd g(u_j \cd i)^3=g(u_j )^{-4}$.
Now we have either $g(v_j )=1$ or $g(v_j \cd i)=1$. Thus at least one of
$g(u_j )^3 =1$ and $g(u_j )^{-4} =1$ must hold. Thus $g(u_j )^{12}=1$
in both cases, which gives $g(a_j )=g(12 u_j )=g(u_j )^{12} =1$.

Since each $x_j$ is the sum of some of the numbers $a_1 \stb a_N$,
it follows that $g(x_j )$ is the product of
some of the numbers $g(a_1 ) \stb g(a_N )$. Thus $g(x_j )=1$ for every
$j=1\stb 4$. However, by $\sum_{j=1}^4 \al _j g(x_j )=0$ this implies
$\sum_{j=1}^4 \al _j =0$, which contradicts the assumption
$\sum_{j=1}^4 \al _j \ne 0$. This contradiction proves that $f\equiv 0$.
$\square$

\bsk
Finally, we present a family of quadrangles depending on a continuous parameter
such that each member of the family has the discrete
Pompeiu property w.r.t. the isometry group.

Let a non-regular triangle $ABC \triangle$ be given in the plane.
The steps of the construction are summarized as follows:

\begin{figure}[h]
\centering
\includegraphics[viewport=0 0 404 234, scale=0.5]{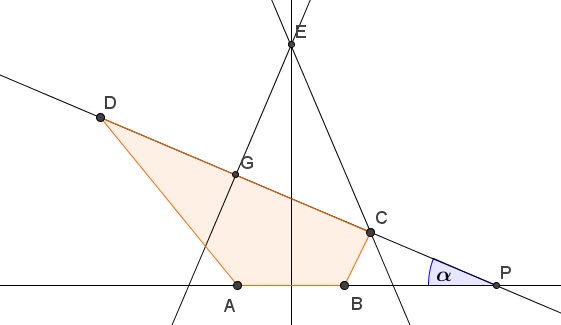}
\caption{A Pompeiu quadrangle belonging to $\alpha=23^{\circ}$.}
\end{figure}

\begin{itemize}
\item since $ABC \triangle$ is non-regular we can suppose that the point $C$ is not on the perpendicular bisector of $AB$; especially $C$ and $B$ are supposed to be on the same side of the perpendicular bisector of $AB$.
\item let $0< \alpha < 45^{\circ}$ be a given angle and choose a point $P$ on the line $AB$ such that $A$ and $P$ are separated by the point $B$ and the angle enclosed by the lines $PB$ and $PC$ is of measure $\alpha$ (see Figure 1).

\item $E$ is the point of the perpendicular bisector of $AB$
such that the line $CE$ intersects the bisector under an angle of measure $\alpha$ (see Figure 1).

\item $G$ is the point on the line $PC$ such that the triangle
$EGC \triangle$ has a right angle at $G$.
Then, necessarily, the perpendicular bisector of $AB$ is the bisector of the
angle of $EGC \triangle$ at the vertex $E$.

\item $D_{\alpha}$ is the reflection of $C$ about the point $G$.
\end{itemize}

\begin{theorem}\label{t4}
The set $H_{\alpha}=\{A, B, C, D_{\alpha}\}$ has the Pompeiu property w.r.t.
the isometry group.
\end{theorem}

\proof
Suppose that the angle $\alpha$ is given,
and let $D:=D_{\alpha}$ for the sake of simplicity.
For any point $P$ let $P'$ be the image of $P$ under the reflection about the perpendicular bisector of $AB$. Then $A'=B$, $B'=A$ and the points $C$, $C'$, $D$ and $D'$ form a symmetric trapezium such that $D'C=CC'=C'D$; see Figure 2. Using that
$$f(A)+f(B)+f(C)+f(D)=0\ \ \textrm{and}\ \ f(A')+f(B')+f(C')+f(D')=0$$
it follows that the alternating sum of the values of $f$ at the vertices of the trapezium $CC'DD'$ vanishes, i.e.
\begin{equation}
\label{symtrap}
f(C)-f(C')+f(D)-f(D')=0.
\end{equation}
Since equation (\ref{symtrap}) holds on any congruent copy of the trapezium $CC'DD'$ we have
\begin{equation}
\label{trick1}
f(C)-f(C')+f(D)-f(D')=0\ \ \textrm{and}\ \ f(C')-f(D)+f(H')-f(C)=0
\end{equation}
as Figure 2 shows: the trapezium $CC'DH'$ comes by a translation $C \mapsto C'$ and a rotation about the point $D$.
\begin{figure}
\centering
\includegraphics[viewport=0 0 493 300, scale=0.5]{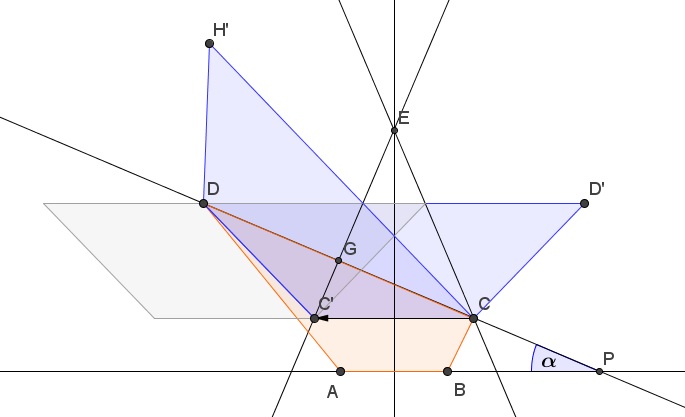}
\caption{The proof of Theorem \ref{t4}.}
\end{figure}
Therefore
\begin{equation}
\label{key}
f(D')=f(H')
\end{equation}
and equation (\ref{key}) holds on any congruent copy of the segment $D'H'$ of measure $r$. This means that $f$ takes the same values at any pair of points having distance $r$. Since any pair of points can be joined by a (finite) chain of circles with radius $r$ it follows that $f$ is a constant function. Especially, the constant must be zero.
$\square$

\begin{small}\noindent
(G. Kiss)\\
{\sc University of Luxembourg, Faculty of Science\\
Mathematical Research Unit\\
E-mail: {\tt gergely.kiss@uni.lu}}

\bsk \noi
(M. Laczkovich)\\
{\sc Department of Analysis, E\"otv\"os Lor\'and University\\
Budapest, P\'azm\'any P\'eter s\'et\'any 1/C, 1117 Hungary\\
E-mail: {\tt laczk@cs.elte.hu}}

\bsk \noi
(Csaba Vincze)\\
{\sc Institute of Mathematics, University of Debrecen\\
P.O.Box 400, Debrecen, 4002 Hungary\\
E-mail: {\tt csvincze@science.unideb.hu}}
\end{small}

\end{document}